\documentclass[12pt]{article}
\usepackage{amssymb,latexsym,amsmath,enumerate,verbatim,amsfonts,amsthm}
\usepackage[active]{srcltx}

\numberwithin{equation}{section}
\newtheorem{theorem}{Theorem}[section]
\newtheorem{definition}{Definition}[section]
\newtheorem{proposition}[theorem]{Proposition}

\newtheorem{lemma}{Lemma}[section]
\newtheorem{corollary}[theorem]{Corollary}

\newtheorem{example}{Example}[section]

\makeatletter
 \def\@biblabel#1{#1.}
\makeatother

\newcommand{\x}{{\bf x}}
\newcommand{\y}{{\bf y}}
\newcommand{\z}{{\bf z}}
\newcommand{\q}{{\bf q}}

\newcommand{\e}{{\bf e}}
\newcommand{\0}{{\bf 0}}

\newcommand{\w}{{\bf w}}

\begin{document}
\title{Properties of Solution Set of Tensor Complementarity Problem}

\author{Yisheng Song\thanks{Corresponding author. School of Mathematics and Information Science and Henan Engineering Laboratory for Big Data Statistical Analysis and Optimal Control, Henan Normal University, XinXiang HeNan, P.R. China, 453007.
Email: songyisheng1@gmail.com. The work was supported by the National Natural Science Foundation of P.R. China (Grant No. 11571905),  Program for Innovative Research Team (in Science and Technology)  in University of Henan Province(14IRTSTHN023).}\quad Gaohang Yu\thanks{School of Mathematics and Computer Sciences, Gannan Normal University, Ganzhou, 341000, China. Email: maghyu@163.com. This author's work was supported by the National Natural Science Foundation of P.R. China (Grant No. 10926029,11001960), NCET Programm of the Ministry of Education (NCET 13-0738), science and technology programm of Jiangxi Education Committee (LDJH12088).}}

\date{}

 \maketitle

%---------------------------------------------------------------------------------Abstract
\begin{abstract}
\noindent  %\vspace{3mm}
     The tensor complementarity problem is a specially  structured nonlinear  complementarity problem,   then it has its particular and nice properties other than ones of the classical nonlinear  complementarity problem.  In this paper, it is proved that  a tensor is an S-tensor if and only if  the tensor complementarity problem is feasible, and each Q-tensor is an S-tensor.  Furthermore, the boundedness of  solution set of  the tensor complementarity problem is equivalent to  the uniqueness of solution for such a problem with zero vector. For the tensor complementarity problem with a strictly semi-positive tensor, we proved the global upper bounds for solution of such a problem. In particular, the upper bounds keep in close contact with the smallest Pareto $H-$($Z-$)eigenvalue. \vspace{3mm}

\noindent {\bf Key words:}\hspace{2mm} Tensor complementarity,  S-tensor,  strictly semi-positive,  feasible vector.  \vspace{3mm}

\noindent {\bf AMS subject classifications (2010):}\hspace{2mm}
47H15, 47H12, 34B10, 47A52, 47J10, 47H09, 15A48, 47H07.
  \vspace{3mm}
\end{abstract}

\section{Introduction}
The nonlinear complementarity problem  was introduced by Cottle in his Ph.D. thesis in 1964. In the last decades, many mathematical workers have concentrated a lot of their energy and attention on this classical  problem because of  a multitude of interesting connections to numerous disciplines and a wide range of important applications in operational research, applied science and technology such as optimization,   economic equilibrium problems,  contact mechanics problems, structural mechanics problem, nonlinear obstacle problems,  traffic equilibrium problems, and discrete-time optimal control. For more detail, see \cite{FP97,HP90,HXQ,FP11} and references therein. Well over a thousand articles and several books have been published on this classical subject, which has developed into a well-established and fruitful discipline in the field of mathematical programming.

The linear complementarity problem is a special case of nonlinear complementarity problem.  It is well-known that the linear complementarity problem has wide and important applications in engineering and economics (Cottle, Pang and Stone \cite{CPS} and Han, Xiu and Qi \cite{HXQ}).  Cottle and Dantzig \cite{CD68} studied the  the existence of solution of  the linear complementarity problem  with the help of the structure of the matrix.   Some relationship between the unique and existence of solution of  the linear complementarity problem and  semi-monotoncity of the matrices were showed by Eaves \cite{E71},  Karamardian \cite{K72}, Pang \cite{P79,P81} and Gowda \cite{S90}, respectively.   Cottle \cite{C80} studied some classes of the complete Q-matrix (a matrix is called the complete Q-matrix iff it and all its principal sub-matrices are Q-matrix), and obtained that each completely Q-matrix is a strictly semi-monotone matrix.%, and obtained the following conclusion about the LCP$(A,\q)$ and the structure of  matrix.

%\begin{Definition} \label{d11}\em A matrix $A$ is called a {\bf Q-matrix} iff the LCP$(A,\q)$ has a solution for all $\q\in\mathbb{R}^n$.  A matrix $A$ is said to be a {\bf completely Q-matrix} iff $A$ and all its principal sub-matrices are Q-matrices.
%\end{Definition}

%\begin{Theorem}\em(Cottle \cite[Theorem 1]{C80})\label{th12} A completely Q-matrix must be a strictly semi-monotone matrix, i.e., for each $\x\geq\0$ and $\x\ne\0$, there exists an index $k\in I_n$ such that $$x_k>0\mbox{ and }\left(A \x\right)_k>0.$$ The reverse was also true.\end{Theorem}

The tensor complementarity problem, as a special type of nonlinear com- plementarity problems, is a new topic emerged from the tensor community, inspired by the growing research on structured tensors. At the same time, the tensor complementarity problem, as a natural extension of the linear complementarity problem seems to have similar  properties to the problem, and to have its particular and nice properties other than ones of the classical linear complementarity problem. So how to obtain the nice properties and their applications of the tensor complementarity problem will be very interesting by means of the  special   structure of  higher order tensor (hypermatrix).

The notion of the tensor complementarity problem is used firstly by Song and Qi \cite{SQ2015}, and they showed the existence of solution for such a problem with some classes of structured tensors. In particular, they showed that the nonnegative tensor complementarity problem has a solution if and only if all principal  diagonal entries of such a tensor are positive. Che, Qi, Wei \cite{CQW} showed the existence of solution for the symmetric positive definite tensor complementarity problem and copositive tensors. Luo, Qi, Xiu \cite{LQX} studied the sparsest solutions to $Z$-tensor complementarity problems. Song and Qi \cite{SQ-2015} studied the solution of the semi-positive tensor complementarity problem, and obtained that a symmetric tensor is (strictly) semi-positive if and only if  it is (strictly) copositive.

In this paper, we will study the properties of solution of the tensor complementarity problem by means of the special structure of tensors.  Namely, it will be proved that the solution set of the tensor complementarity problem is bounded if and only if such a tensor is a R$_0$-tensor. We will present the global upper bounds for solution of the tensor complementarity problem with a strictly semi-positive tensor.   % Moreover, an example is given to verify that the  function  $F(x)=\q + \mathcal{A}\x^{m-1}$ is strictly semi-positive, not pseudo-monotone.  %So our result is clearly different from Theorem \ref{th11}.

In section 2, we will give some notions, definitions and basic facts, which will be used in the context.  In section 3, we will introduce the concept of S-tensor and give its two equivalent definitions by means of the tensor complementarity problem.  Subsequently, it is proved that each Q-tensor is an S-tensor, and a tensor is an S-tensor if and only if  the corresponding tensor complementarity problem is feasible.  Furthermore, the solution set of  the tensor complementarity problem is bounded if and only if  such a problem with zero vector has a unique solution.  Finally, with the help of the smallest Pareto $H-$($Z-$)eigenvalue, we will discuss the global upper bounds for solution of the tensor complementarity problem with a strictly semi-positive tensor.	
%a real tensor $\mathcal{A}$ is completely S-tensor if and only if it is a  semi-positive tensor.   That is,  the following four classes of tensors are equivalent:
%$$\begin{aligned}\mbox{completely}&  \mbox{ S-tensor }\Leftrightarrow\  \mbox{ completely}& &\mbox{Q-tensor }\\
%\Updownarrow& &\Updownarrow&\\
% \mbox{ completely}&  \mbox{ R-tensor }\Leftrightarrow\ \ \ \ \mbox{  strictly}& &\mbox{semi-positive tensor}
% \end{aligned}$$

\section{Preliminaries}
%\hspace{4mm}

Throughout this paper, we use small letters $x, y, v, \alpha, \cdots$, for scalars, small bold
letters $\x, \y,  \cdots$, for vectors, capital letters $A, B,
\cdots$, for matrices, calligraphic letters $\mathcal{A}, \mathcal{B}, \cdots$, for
tensors.  All the tensors discussed in this paper are real.  Let $I_n := \{ 1,2, \cdots, n \}$, and $\mathbb{R}^n:=\{(x_1, x_2,\cdots, x_n)^\top;x_i\in  \mathbb{R},  i\in I_n\}$, $\mathbb{R}^n_{+}:=\{x\in \mathbb{R}^n;x\geq\0\}$, $\mathbb{R}^n_{-}:=\{\x\in \mathbb{R}^n;x\leq\0\}$, $\mathbb{R}^n_{++}:=\{\x\in \mathbb{R}^n;x>\0\}$, where $\mathbb{R}$ is the set of real numbers, $\x^\top$ is the transposition of a vector $\x$, and $\x\geq\0$ ($\x>\0$) means $x_i\geq0$ ($x_i>0$) for all $i\in I_n$.

Let $F$ be a nonlinear function from $\mathbb{R}^n$ into itself.  The nonlinear complementarity problem is to  find a vector $ \x \in \mathbb{R}^n$ such that
$$ \x \geq \0, F(\x) \geq \0, \mbox{ and }\x^\top F(\x) = 0, $$
or to show that no such vector exists.
Let $A = (a_{ij})$ be an $n \times n $ real matrix and $F(\x)=\q + A\x$.  Then such a nonlinear complementarity problem  is called the linear complementarity problem, denoted by LCP $(A, \q)$, i.e.,   to find $ \x \in \mathbb{R}^n$ such that
	$$	\x \geq \0, \q + A\x \geq \0, \mbox{ and }\x^\top (\q + A\x) = 0,$$or to show that no such vector exists.

In 2005, Qi \cite{Qi} introduced the concept of
positive (semi-)definite symmetric tensors.   A real $m$th order $n$-dimensional tensor (hypermatrix) $\mathcal{A} = (a_{i_1\cdots i_m})$ is a multi-array of real entries $a_{i_1\cdots
	i_m}$, where $i_j \in I_n$ for $j \in I_m$. Denote the set of all
real $m$th order $n$-dimensional tensors by $T_{m, n}$. Then $T_{m,
	n}$ is a linear space of dimension $n^m$. Let $\mathcal{A} = (a_{i_1\cdots
	i_m}) \in T_{m, n}$. If the entries $a_{i_1\cdots i_m}$ are
invariant under any permutation of their indices, then $\mathcal{A}$ is
called a {\bf symmetric tensor}.  Denote the set of all real $m$th
order $n$-dimensional symmetric tensors by $S_{m, n}$. Then $S_{m, n}$ is a
linear subspace of $T_{m, n}$.  We
denote the zero tensor in $T_{m, n}$ by $\mathcal{O}$.   Let $\mathcal{A} =(a_{i_1\cdots i_m}) \in
T_{m, n}$ and $\x \in \mathbb{R}^n$. Then $\mathcal{A} \x^{m-1}$ is a vector in $\mathbb{R}^n$ with
its $i$th component as
$$\left(\mathcal{A} \x^{m-1}\right)_i: = \sum_{i_2, \cdots, i_m=1}^n a_{ii_2\cdots
	i_m}x_{i_2}\cdots x_{i_m}$$ for $i \in I_n$.
Then $\mathcal{A} \x^m$ is a homogeneous
polynomial of degree $m$, defined by
$$\mathcal{A} \x^m:= \x^\top(\mathcal{A} \x^{m-1})= \sum_{i_1,\cdots, i_m=1}^n a_{i_1\cdots i_m}x_{i_1}\cdots
x_{i_m}.$$
A tensor $\mathcal{A} \in T_{m, n}$ is called {\bf positive
	semi-definite} if for any vector $\x \in \mathbb{R}^n$, $\mathcal{A} \x^m \ge 0$,
and is called {\bf positive definite} if for any nonzero vector $\x
\in \mathbb{R}^n$ and an even number $m$, $\mathcal{A} \x^m > 0$.
Recently, miscellaneous structured tensors are widely studied, for example,
Zhang, Qi and Zhou \cite{ZQZ} and
Ding, Qi and Wei \cite{DQW} for M-tensors, Song and Qi \cite{SQ-15} for P-(P$_0$)tensors and B-(B$_0$)tensors, Qi and Song \cite{QS} for  B-(B$_0$)tensors,  Song and Qi \cite{SQ1} for infinite and finite dimensional Hilbert tensors, Song and Qi \cite{SQ} for E-eigenvalues of weakly symmetric nonnegative tensors and references therein.

\begin{definition} \label{d22}  Let $\mathcal{A}  = (a_{i_1\cdots i_m}) \in T_{m, n}$.   \begin{itemize}
		\item[(i)] The {\bf tensor complementarity problem}, denoted by TCP $(\mathcal{A},\q)$, is to find  $ \x \in \mathbb{R}^n$ such that
		\begin{align}\label{eq:21} \x \geq& \0\\
		\label{eq:22} \q + \mathcal{A}\x^{m-1} \geq& \0\\
		\label{eq:23}\x^\top (\q + \mathcal{A}\x^{m-1}) =& 0
		\end{align}or to show that no such vector exists.
		\item[(ii)] A vector $\x$ is said to be
		{\bf (strictly)  feasible } iff $\x$  satisfies the  inequality (\ref{eq:21}) and (strict) inequality (\ref{eq:22}).	
		\item[(iii)] The TCP $(\mathcal{A}, \q)$ is said to be {\bf (strictly) feasible} iff a (strictly)  feasible vector exists.
		\item[(iv)] The set of all feasible vector of the TCP $(\mathcal{A}, \q)$ is said to be its {\bf feasible region}.
		\item[(v)] The TCP $(\mathcal{A}, \q)$ is said to be {\bf solvable} iff there is a feasible vector  satisfying the equation (\ref{eq:23}).
		\item[(vi)] $\mathcal{A}$ is called a {\bf Q-tensor} (\cite{SQ2015}) iff the TCP $(\mathcal{A}, \q)$ is  solvable for all $\q\in\mathbb{R}^n$.
		\item[(vii)] $\mathcal{A}$ is called a {\bf R$_0$-tensor} (\cite{SQ2015}) iff the TCP $(\mathcal{A}, \0)$ has unique solution.
		\item[(viii)] $\mathcal{A}$ is called a {\bf R-tensor} (\cite{SQ2015}) iff it is a R$_0$-tensor and  the TCP $(\mathcal{A}, \q)$ has unique solution for $\q=(1,1,\cdots,1)^\top$.
	\end{itemize}
\end{definition}
Let $\w=\q + \mathcal{A}\x^{m-1}$.  Then a feasible vector $\x$ of the TCP $(\mathcal{A}, \q)$ is its solution  if and only if \begin{equation}\label{eq:25} x_iw_i=x_iq_i+x_i\left(\mathcal{A}\x^{m-1}\right)_i=0\mbox{ for all }i\in I_n.\end{equation}

Recently, Song and Qi \cite{SQ2015} extended the concepts of  (strictly) semi-monotone matrices  to  (strictly) semi-positive tensors.

\begin{definition} \label{d23}
	Let $\mathcal{A}  = (a_{i_1\cdots i_m}) \in T_{m, n}$.   $\mathcal{A}$ is said to be\begin{itemize}
		\item[(i)] {\bf semi-positive} iff for each $\x\geq0$ and $\x\ne\0$, there exists an index $k\in I_n$ such that $$x_k>0\mbox{ and }\left(\mathcal{A} \x^{m-1}\right)_k\geq0;$$
		\item[(ii)]  {\bf strictly semi-positive} iff for each $\x\geq\0$ and $\x\ne\0$, there exists an index $k\in I_n$ such that $$x_k>0\mbox{ and }\left(\mathcal{A} \x^{m-1}\right)_k>0;$$
		\item[(iii)] a {\bf P-tensor}(Song and Qi \cite{SQ-15}) iff for each $\x$ in $\mathbb{R}^n$ and $\x\ne\0$, there
		exists $i \in I_n$ such that
		$$x_i \left(\mathcal{A} \x^{m-1}\right)_i > 0;$$
		\item[(iv)] a {\bf P$_0$-tensor}(Song and Qi \cite{SQ-15}) iff for every $\x$ in $\mathbb{R}^n$ and $\x\ne\0$,, there
		exists $i \in I_n$ such that $x_i \not = 0$ and
		$$x_i \left(\mathcal{A} \x^{m-1}\right)_i \geq 0;$$
      \item[(v)] {\bf copositive }(Qi \cite{LQ5}) if $\mathcal{A}\x^m\geq0$ for all $\x\in \mathbb{R}^n_+$;
		\item[(vi)] {\bf strictly copositive}(Qi \cite{LQ5}) if  $\mathcal{A}\x^m>0$ for all $\x\in \mathbb{R}^n_+\setminus\{\0\}$.
	\end{itemize}
\end{definition}

\begin{lemma}\label{le21} (Song and Qi \cite[Corollary 3.3, Theorem 3.4]{SQ2015}) Each strictly semi-positive tensor must be a R-tensor, and each R-tensor must be a Q-tensor. A semi-positive R$_0$-tensor is a Q-tensor.
\end{lemma}

Song and Qi \cite{SQ2013} introduced the concept of Pareto $H-$($Z-$)eigenvalue and used it to portray the (strictly) copositive tensor. The number and computation of Pareto $H-$($Z-$)eigenvalue see Ling, He and Qi \cite{LHQ}, Chen, Yang and Ye \cite{CYY}.

\begin{definition} \label{d24}
	Let $\mathcal{A}  = (a_{i_1\cdots i_m}) \in T_{m, n}$.   A real number $\mu$ is said to be
\begin{itemize}
\item[(i)]{\bf Pareto $H-$eigenvalue} of $\mathcal{A}$ iff there is a non-zero vector $\x\in \mathbb{R}^n$ satisfying \begin{equation}\label{eq:26} \mathcal{A}\x^m=\mu \x^T\x^{[m-1]}, \
\mathcal{A}\x^{m-1}-\mu \x^{[m-1]} \geq0,\ \x\geq 0.\end{equation}
\item[(ii)] {\bf Pareto $Z-$eigenvalue} of $\mathcal{A}$ iff there is a non-zero vector $\x\in \mathbb{R}^n$ satisfying \begin{equation}\label{eq:27}  \mathcal{A}\x^m=\mu (\x^T\x)^{\frac{m}2}, \
\mathcal{A}\x^{m-1}-\mu (\x^T\x)^{\frac{m}2-1}\x \geq0,\  \x\geq 0. \end{equation}
	\end{itemize}
\end{definition}

\begin{lemma}\label{le22} (Song and Qi \cite[Theorem 3.1,3.3, Corollary 3.5]{SQ2013}) Let $\mathcal{A}\in S_{m,n}$ be strictly copositive. Then
\begin{itemize}
\item[(i)] $\mathcal{A}$ has at least one Pareto $H$-eigenvalue $\lambda(\mathcal{A}):=\min\limits_{\x\geq0 \atop \|\x\|_m=1 }\mathcal{A}\x^m$ and
\begin{equation}\label{eq:28}\lambda(\mathcal{A})=\min\{\lambda; \lambda \mbox{ is Pareto $H$-eigenvalue of  }\mathcal{A}\}>0;
\end{equation}
\item[(ii)] $\mathcal{A}$ has at least one Pareto $Z$-eigenvalue $\mu(\mathcal{A}):=\min\limits_{\x\geq0 \atop \|\x\|_2=1 }\mathcal{A}\x^m$ and
\begin{equation}\label{eq:29}\mu(\mathcal{A})=\min\{\mu; \mu \mbox{ is Pareto $Z$-eigenvalue of  }\mathcal{A}\}>0.
\end{equation}\end{itemize}
\end{lemma}

%\begin{Definition} \label{d24}\em Let $\mathcal{A}  = (a_{i_1\cdots i_m}) \in T_{m, n}$. In homogeneous polynomial $\mathcal{A}\x^m$, if we let some (but not all) $x_i$ be zero, then we have a less variable homogeneous polynomial, which defines a lower dimensional tensor. We call such a lower dimensional tensor a {\bf principal sub-tensor} of $\mathcal{A}$.
%We denote the principal subtensor of   $\mathcal{A}$ by $\mathcal{A}^{|N|}$, where  $\emptyset\ne N\subset I_n=\{1, 2,  \cdots, n\}$ and $|N|$ is  the cardinality of $N$.  That is,
%$$\mathcal{A}^{|N|} = (a_{i_1\cdots i_m}),\mbox{ for all } i_1, i_2, \cdots, i_m\in \mathcal{N}.$$
%\end{Definition}
% Clearly, $\mathcal{A}^{|N|}$ is a tensor of order $m$ and dimension $|N|$ and the principal sub-tensor $\mathcal{A}^{|N|}$  is just $\mathcal{A}$ itself when $N=I_n$.
%The concept were first introduced and used by Qi \cite{Qi} to the higher order symmetric tensor.

% \begin{Definition} \label{d25}\em	Let $Y$ denote a fixed class of higher order tensor.   A tensor $\mathcal{A}$ is said to be {\bf completely $Y$} iff $\mathcal{A}$ and all its principal sub-tensor belong to $Y$.\end{Definition}

%The following conclusions are  obvious.

%\begin{Proposition} \em Each  semi-positive tensor is completely  semi-positive tensor. Each strictly semi-positive tensor is completely strictly semi-positive tensor.\end{Proposition}

\section{Solution of the TCP$(\mathcal{A}, \q)$}
%\hspace{4mm}
\subsection{\bf S-tensor and feasible solution of the TCP$(\mathcal{A}, \q)$ }

We first introduce the concept of the S-tensor, which is a natural extension of S-matrix \cite{CPS}.
\begin{definition} \label{d21}
	Let $\mathcal{A}  = (a_{i_1\cdots i_m}) \in T_{m, n}$.   $\mathcal{A}$ is said to be a\begin{itemize}
		\item[(i)] {\bf S-tensor} iff the system $$\mathcal{A}\x^{m-1}>\0, \ \x>\0$$ has a solution;
		\item[(ii)]  {\bf S$_0$-tensor} iff the system $$\mathcal{A}\x^{m-1}\geq\0, \ \x\geq\0,\ \x\ne\0$$ has a solution.\end{itemize}
\end{definition}

\begin{proposition}  \label{pro:21}
	Let $\mathcal{A}\in T_{m,n}$. Then $\mathcal{A}$ is a  S-tensor if and only if the system \begin{equation}\label{eq:2 1}\mathcal{A}\x^{m-1}>\0, \ \x\geq\0\end{equation} has a solution.
	%\item[(ii)]  $\mathcal{A}$ is a  S$_0$-tensor if and only if the system \begin{equation}\label{eq:2-1}\mathcal{A}\x^{m-1}\geq\0, \ \x>\0\end{equation} has a solution.\end{itemize}
\end{proposition}
{\it Proof}
	It follows from Definition \ref{d21} that the necessity is obvious.  Now we show the  sufficiency. In fact, if there exists $\y$ such that $$\mathcal{A}\y^{m-1}>\0, \ \y\geq\0,$$
	clearly, $\y\ne \0$. Since $\mathcal{A}\y^{m-1}$ is continuous on $\y$,  it follows that $\mathcal{A}(\y+t\e)^{m-1}>0$ for some small enough $t>0$, where $\e=(1,1,\cdots,1)^\top$. It is obvious that  $\y+t\e>0$. So $\mathcal{A}$ is an S-tensor.% from the locally sign-preserving property of continuous function (ii) The sufficiency is obvious. The proof of the necessity is the same as 	ones of the  sufficiency of (i) with appropriate changes in the inequalities, which is a repetitive work, we omit it.
\qed

Now, by means of the solution of the TCP $(\mathcal{A},\q)$, we may give the following equivalent definition of S-tensor.

\begin{theorem}  \label{th:22}
	Let $\mathcal{A}\in T_{m,n}$.  Then $\mathcal{A}$  is a  S-tensor if and only if the TCP $(\mathcal{A}, \q)$ is  feasible for all $\q\in\mathbb{R}^n$. Meanwhile, each Q-tensor must be an S-tensor.
\end{theorem}
{\it Proof} Let $\mathcal{A}$ is a  S-tensor. Then it follows from Definition \ref{d21} of S-tensor that  there exists $\y$ such that $$\mathcal{A}\y^{m-1}>\0, \ \y>\0.$$
	So for each $\q\in\mathbb{R}^n$, there is $t>0$ such  that
	$$\mathcal{A}(\sqrt[m-1]{t}\y)^{m-1}=t\mathcal{A}\y^{m-1}\geq-\q.$$
	Clearly, $\sqrt[m-1]{t}\y> \0$, and so $\sqrt[m-1]{t}\y$ is a feasible vector of the TCP $(\mathcal{A},\q)$.
	
	On the other hand,  if the TCP $(\mathcal{A}, \q)$ is  feasible for all $\q\in\mathbb{R}^n$,  we take $\q<\0$.  Let $\z$ is a feasible solution  of the TCP $(\mathcal{A},\q)$. Then
	$$\z\geq\0\mbox{ and }\q + \mathcal{A}\z^{m-1} \geq \0,$$
	and hence $$ \mathcal{A}\z^{m-1} \geq -\q > \0.$$
	So $\z$ is a solution of the system (\ref{eq:2 1}). It follows from Proposition \ref{pro:21} that  $\mathcal{A}$ is an S-tensor.\qed		

\subsection{\bf R$_0$-tensor and boundedness of solution set of the TCP $(\mathcal{A},\q)$  }

%Now we give the boundedness of solution set of the TCP $(\mathcal{A},\q)$.

\begin{theorem}  \label{th:33}
	Let $\mathcal{A}\in T_{m,n}$.  Then the following three conclusions are equivalent:
	\begin{itemize}
		\item[(i)]  $\mathcal{A}$ is R$_0$-tensor;
		\item[(ii)]  For each $\q\in \mathbb{R}^n$ and each $t, s\in\mathbb{R}$ with $t>0$, the set  $$\Gamma(\q,s,t)=\{\x\geq\0; \q+\mathcal{A}\x^{m-1}\geq\0\mbox{ and }\x^\top\q+t\mathcal{A}\x^m\leq s\}$$ is bounded;
		\item[(iii)] For each $\q\in \mathbb{R}^n$,  the solution set of the TCP $(\mathcal{A},\q)$ is bounded.
	\end{itemize}
\end{theorem}
{\it Proof} (i) $\Rightarrow$ (ii).  Suppose that there exist $\q'\in \mathbb{R}^n$, $s'\in\mathbb{R}$ and $t'>0$ such that the set $\Gamma(\q',s',t')$ is not bounded. Let a sequence $\{\x^k\}\subset\Gamma(\q',s',t')$ be an unbounded sequence. Then the sequence $\{\frac{\x^k}{\|\x^k\|}\}$ is bounded, and so there exists $\x'\in \mathbb{R}^n$ and a subsequence $\{\frac{\x^{k_j}}{\|\x^{k_j}\|}\}$ such that $$\lim_{j\to\infty}\frac{\x^{k_j}}{\|\x^{k_j}\|}=\x'\ne\0 \mbox{ and } \lim_{j\to\infty}\|\x^{k_j}\|=\infty.$$ From the definition of $\Gamma(\q',s',t')$, it follows that
	\begin{equation}\label{eq:32}\frac{\q'}{\|\x^{k_j}\|^{m-1}}+\mathcal{A}(\frac{\x^{k_j}}{\|\x^{k_j}\|})^{m-1}\geq\0\mbox{ and }\frac{\x^\top\q'}{\|\x^{k_j}\|^m}+t'\mathcal{A}(\frac{\x^{k_j}}{\|\x^{k_j}\|})^m\leq \frac{s'}{\|\x^{k_j}\|^m},\end{equation}
	and hence, by the continuity of $\mathcal{A}\x^m$ and $\mathcal{A}\x^{m-1}$, let $j\to\infty$,
	$$\mathcal{A}(\x')^{m-1}\geq\0\mbox{ and }\mathcal{A}(\x')^m\leq0.$$
	Since $\x'\geq\0$,  we have $$\mathcal{A}(\x')^m=(\x')^\top\mathcal{A}(\x')^{m-1}\geq\0.$$
	Thus, $\mathcal{A}(\x')^m=0$, and hence, $\x'$ is a nonzero solution of the TCP$(\mathcal{A},\0)$. This contradicts the assumption that $\mathcal{A}$ is R$_0$-tensor.
	
	(ii) $\Rightarrow$ (iii).   It follows from the definition of $\Gamma(\q,s,t)$ that
	$$\Gamma(\q,0,1)=\{\x\geq\0; \q+\mathcal{A}\x^{m-1}\geq\0\mbox{ and }\x^\top(\q+\mathcal{A}\x^{m-1})= 0\}.$$
	That is, $\Gamma(\q,0,1)$ is the solution set of the TCP$(\mathcal{A},\q)$. So the conclusion follows.
	
	(iii) $\Rightarrow$ (i). Suppose $\mathcal{A}$ is not R$_0$-tensor. Then the TCP$(\mathcal{A},\0)$ has a nonzero solution $\x^*$, and so $\x^*\in \Gamma(\0,0,1)$.
	Since $\x^*\ne\0$, $\tau\x^*\in \Gamma(\0,0,1)$ for all $\tau>0$. Therefore, the set $\Gamma(\0,0,1)$ is not bounded. This contradicts the assumption (iii). So $\mathcal{A}$ is R$_0$-tensor.\qed

It is known that each strictly semi-positive tensor is an R-tensor and each P-tensor is a strictly semi-positive tensor (Song and Qi \cite{SQ2015}). The following conclusions are obvious.

\begin{corollary}  \label{co:33}
	Let $\mathcal{A}$ be a strictly semi-positive tensor.  Then for each $\q\in \mathbb{R}^n$,  the solution set of the TCP $(\mathcal{A},\q)$ is bounded.
\end{corollary}

\begin{corollary}  \label{co:33}
	Let $\mathcal{A}$ be a  P-tensor.  Then for each $\q\in \mathbb{R}^n$,  the solution set of the TCP$(\mathcal{A},\q)$ is bounded.
\end{corollary}

\subsection{\bf Solution of TCP$(\mathcal{A},\q)$ with strictly semi-positive tensors}

In this section, we discuss the global upper  bound for solution of TCP$(\mathcal{A},\q)$ with strictly semi-positive and symmetric tensor $\mathcal{A}$. Song and Qi \cite{SQ-2015} showed the following conclusion about a symmetric tensor.

\begin{lemma}\label{le31} (Song and Qi \cite[Theorem 3.2, 3.4]{SQ-2015}) Let $\mathcal{A}  = (a_{i_1\cdots i_m}) \in T_{m,n}$. Then a symmetric tensor $\mathcal{A}$ is strictly semi-positive if and only if  it is strictly copositive. Moreover, the TCP $(\q, \mathcal{A})$ has a unique solution $\0$ for  $\q \geq \0$ when $\mathcal{A}$ is strictly semi-positive.
\end{lemma}

\begin{theorem}\label{th:34} Let $\mathcal{A}  = (a_{i_1\cdots i_m}) \in S_{m,n}$ be strictly semi-positive. If $\x$ is a solution of the TCP $(\q, \mathcal{A})$,  then
\begin{equation}\label{eq:33}
\|\x\|_m^{m-1}\leq\frac{\|(-\q)_+\|_\frac{m}{m-1}}{\lambda(\mathcal{A})},
\end{equation}
where $\lambda(\mathcal{A})$ is defined  in Lemma \ref{le22} (i), $\x_+:=\left(\max\{x_1,0\},\max\{x_2,0\},\cdots,\max\{x_n,0\}\right)^\top$ and $\|\x\|_m:=\left(\sum\limits_{i=1}^n|x_i|^m\right)^{\frac1m}$.
\end{theorem}
{\it Proof} It follows from Lemma \ref{le22} and \ref{le31} that
\begin{equation}\label{eq:34}\lambda(\mathcal{A})=\min\{\lambda; \lambda \mbox{ is Pareto $H$-eigenvalue of  }\mathcal{A}\}=\min\limits_{\y\geq0 \atop \|\y\|_m=1 }\mathcal{A}\y^m>0.
\end{equation}
Since $\x$ is a solution of the TCP $(\q, \mathcal{A})$, we have
$$\mathcal{A}\x^{m}-\x^\top(-\q)=\x^\top(\mathcal{A}\x^{m-1}+\q)=0,\ \mathcal{A}\x^{m-1}+\q\geq\0\mbox{ and }\x\geq\0.$$ Suppose that $\q\geq\0.$ Then $\x=\0$ by Lemma \ref{le31}, the conclusion is obvious. So we may assume that $\q$ is not nonnegative, then $\x\ne\0$ (suppose not, $\x=\0$,   $\mathcal{A}\x^{m-1}+\q=\q$, which contradict to the fact that $\mathcal{A}\x^{m-1}+\q\geq\0$). Therefore, we have $$\frac{\mathcal{A}\x^m}{\|\x\|^m_m}=\mathcal{A}\left(\frac{\x}{\|\x\|_m}\right)^m\geq\lambda(\mathcal{A})>0.$$
 Thus, we have
$$0<\|\x\|^m_m\lambda(\mathcal{A})\leq\mathcal{A}\x^{m}=\x^\top(-\q)\leq \x^\top(-\q)_+\leq\|\x\|_m\|(-\q)_+\|_\frac{m}{m-1}.$$
The desired conclusion follows. \qed

\begin{theorem}\label{th:35} Let $\mathcal{A}  = (a_{i_1\cdots i_m}) \in S_{m,n}$ be strictly semi-positive. If $\x$ is a solution of the TCP $(\q, \mathcal{A})$,  then
\begin{equation}\label{eq:35}
\|\x\|_2^{m-1}\leq\frac{\|(-\q)_+\|_2}{\mu(\mathcal{A})},
\end{equation}
where $\mu(\mathcal{A})$ is defined by (\ref{eq:29}) in Lemma \ref{le22}.
\end{theorem}
{\it Proof} It follows from Lemma \ref{le22} and \ref{le31} that
\begin{equation}\label{eq:36}\mu(\mathcal{A})=\min\{\mu; \mu \mbox{ is Pareto $Z$-eigenvalue of  }\mathcal{A}\}=\min\limits_{\y\geq0 \atop \|\y\|_2=1 }\mathcal{A}\y^m>0.
\end{equation}
Similarly, we also  may assume that $\q$ is not nonnegative, then $\x\ne\0$, and so $$\frac{\mathcal{A}\x^m}{\|\x\|^m_2}=\mathcal{A}\left(\frac{\x}{\|\x\|_2}\right)^m\geq\mu(\mathcal{A})>0.$$
 Thus, we have
$$0<\|\x\|^m_2\mu(\mathcal{A})\leq\mathcal{A}\x^{m}=\x^\top(-\q)\leq\x^\top(-\q)_+\leq\|\x\|_2\|(-\q)_+\|_2.$$
The desired conclusion follows. \qed

We now introduce a quantity for a  strictly semi-positive tensor $\mathcal{A}$.
\begin{equation} \label{beta}
\beta(\mathcal{A}):=\min_{\x\geq0 \atop \x \|_\infty=1}\max_{i \in I_n}x_i(\mathcal{A}\x^{m-1})_i,
\end{equation}
where $\|\x\|_\infty:=\max\{|x_i|;i\in I_n\}$.
It follows from the definition of strictly semi-positive tensor that $\beta(\mathcal{A})>0$. Then the following equation is well defined in Theorem \ref{th:36}.

\begin{theorem}\label{th:36} Let $\mathcal{A}  = (a_{i_1\cdots i_m}) \in T_{m,n}$ be strictly semi-positive. If $\x$ is a solution of the TCP $(\q, \mathcal{A})$,  then
\begin{equation}\label{eq:37}
\|\x\|_\infty^{m-1}\leq\frac{\|(-\q)_+\|_\infty}{\beta(\mathcal{A})}.
\end{equation}
\end{theorem}
{\it Proof} Suppose that $\q\geq\0.$ Then $\x=\0$ by Lemma \ref{le31}, the conclusion is obvious. So we may assume that $\q$ is not nonnegative, similarly to the proof technique of Theorem \ref{th:34}, we have $\0\leq\x\ne\0$. It follows from the definition of $\beta(\mathcal{A})$ and (\ref{eq:25}) that
\begin{align}0<\|\x\|^m_\infty\beta(\mathcal{A})\leq &\max_{i \in I_n}x_i(\mathcal{A}\x^{m-1})_i\nonumber
 = \max_{i \in I_n}x_i(-q)_i\\\leq &\max_{i \in I_n}x_i((-q)_+)_i\leq\|x\|_\infty\|(-q)_+\|_\infty.\nonumber
\end{align}
The desired conclusion follows. \qed

It is well known that if a nonlinear complementarity problem with pseudo-monotone and continuous function $F$ has a strictly feasible point $x^*$ (i.e., $\x^*\geq0,\ F(\x^*)>0$), then it has a solution (\cite[Theorem 2.3.11]{HXQ}). A function $F$ from $\mathbb{R}^n_+$ into itself  is called pseudo-monotone iff for all vectors $\x,\y \in \mathbb{R}^n_+$,
$$( \x -\y )^\top F(\y) \geq 0\ \Rightarrow\ ( \x -\y )^\top F(\x)\geq0.$$
Now we give an example to certify the function $F$ deduced by a strictly semi-positive tensor is not pseudo-monotone. However,  the corresponding nonlinear complementarity problem has a solution by Lemma \ref{le21}.

\begin{example}\label{e21}
	Let $\mathcal{A}  = (a_{i_1i_2i_3}) \in T_{3, 2}$ and $a_{111}=1$, $a_{122}=1$, $a_{211}=1$, $a_{221}=-2$, $a_{222}=1$ and all other $a_{i_1i_2i_3}=0$. Then
	$$\mathcal{A} \x^2=\left(\begin{aligned}x_1^2&+x_2^2\\
	x_1^2&-2x_1x_2+x_2^2\end{aligned}\right).$$ Clearly,  $\mathcal{A}$ is strictly semi-positive, and so it is a Q-tensor.
	
	Let $F(\x)=\mathcal{A}\x^{2}+\q$ for $\q=(-\frac32,-\frac12)^\top$.  Then $F$ is not pseudo-monotone. In fact,  $$F(\x)=\mathcal{A} \x^2+\q=\left(\begin{aligned}x_1^2&+x_2^2-\frac32\\
	(x_1^2&-x_2)^2-\frac12\end{aligned}\right).$$  Take $\x=(1,0)^\top$ and $\y=(1,1)^\top$. Then $$\x-\y=\left(\begin{aligned}0\\
	-1\end{aligned}\right),\ F(\x)=\left(\begin{aligned}-\frac12\\
	\frac12\end{aligned}\right)\mbox{ and }F(\y)=\left(\begin{aligned}\frac12\\
	-\frac12\end{aligned}\right).$$ Clearly, we have $$(\x-\y)^\top F(\y)=0\times\frac12-1 \times(- \frac12)>0.$$  However, $$(\x-\y)^\top F(\x)=0\times(-\frac12)-1\times \frac12<0,$$ and hence $F$ is not pseudo-monotone.
\end{example}

%%%%%%%%%%%%%%%%%%%%%%%%%%%%%%%%%%%%%%%%%%%%%%%%%%%%%%%%%%%%%%%%%%%%%%%%%%%%%%%%%%%%%%%%%%%%%%%%%%%%%%%%%%%%%%%%%%%%%%%%%%%%%%%%%%%%%%%%%%%%%%%%%%%%%%%%%%%%
%%%%%%%%%%%%%%%%%%%%%%%%%%%%%%%%%%%%%%%%%%%%%%%%%%%%%%%%%%%%%%%%%%%%%%%%%%%%%%%%%%%%%%%%%%%%
\section{Perspectives}
There are more research topics on solution of the tensor complementarity problem for further research.

 We show the global upper bounds for solution of TCP$(\mathcal{A},\q)$ with strictly semi-positive (symmetric) tensor $\mathcal{A}$.
Then we have the following questions for further research. \begin{itemize}
  \item[(i)]  Do it has a positive lower bound?
  \item[(ii)] Are the above upper bounds is the smallest?
  \item[(iii)] For some other structured tensors such as Z-tensors, B-tensors, H-tensors, M-tensors and so on, whether do they have similar upper bounds or not?
  \end{itemize}

    We define a quantity $\beta(\mathcal{A})$ for a  strictly semi-positive tensor $\mathcal{A}$. Then does the quantity have a  global (or local) upper bound?

%%%%%%%%%%%%%%%%%%%%%%%%%%%%%%%%%%%%%%%%%%%%%%%%%%%%%%%%%%%%%%%%

\section{Conclusions}

In this paper,  We discuss the equivalent relationships between the feasible solution of the tensor complementarity problem and S-tensor, the boundedness of solution of the tensor complementarity problem and R$_0$-tensor. By means of two constants,  keep in close contact with the smallest Pareto $H-$($Z-$)eigenvalue, we proved the global upper bounds for solution of the tensor complementarity problem with a strictly semi-positive and symmetric tensor. A new quantity  is introduced for a strictly semi-positive tensor. For a strictly semi-positive tensor (not symmetric),  the global upper bounds for such a tensor complementarity problem is proved by means of such a new quantity.

%%%%%%%%%%%%%%%%%%%%%%%%%%%%%%%%%%%%%%%%%%%%%%%%%%%%%%%%%%%%%%%%

\section*{\bf Acknowledgment}
%\begin{acknowledgements}
	The authors would like to thank the editors and anonymous
	referees for their valuable suggestions which helped us to improve this manuscript. %The first author's work was supported by the National Natural Science Foundation of P.R. China (Grant No. 11571095),  Program for Innovative Research Team (in Science and Technology)  in University of Henan Province(14IRTSTHN023). The second author's work was supported by the National Natural Science Foundation of P.R. China (Grant No. 10926029,11001960), NCET Programm of the Ministry of Education (NCET 13-0738), science and technology programm of Jiangxi Education Committee (LDJH12088).
%\end{acknowledgements}

%\section*{\bf Acknowledgment}
%The authors would like to thank the anonymous
%referees for their valuable suggestions which helped us to improve this manuscript.

%\bibliographystyle{amsplain}

\end{document}